\documentclass[11pt]{article}

\usepackage{epic,latexsym,amssymb}
\usepackage{color}
\usepackage{tikz}
\usepackage{amsfonts,epsf,amsmath}
\usepackage{epic,latexsym,amssymb,xcolor}
\usepackage{xcolor}
\usepackage{tikz}
\usepackage{amsfonts,epsf,amsmath,leftidx}
\usepackage{float}
\usepackage{mathrsfs}
\usetikzlibrary{arrows}

\usepackage[bottom=1.5in, top=1.2in, left=1.5in, right=1.5in]{geometry}

\newtheorem{thm}{Theorem}
\newtheorem{lem}{Lemma}

\newtheorem{claim}{Claim}

\newtheorem{cor}{Corollary}

\newtheorem{conj}{Conjecture}

\newcommand{\qed}{$\Box$}

\newcommand{\smallqed}{{\tiny ($\Box$)}}

\newcommand{\cF}{{\cal F}}

\newcommand{\proof}{\noindent\textbf{Proof. }}

\let\oldenumerate\enumerate
\renewcommand{\enumerate}{
  \oldenumerate
  \setlength{\itemsep}{0pt}
  \setlength{\parskip}{0pt}
  \setlength{\parsep}{0pt}
}

\begin{document}

\title{Partial Domination and Irredundance Numbers in Graphs}

\author{Pawaton Kaemawichanurat$^{\ddagger}$ and Odile Favaron
\\ \\
$^{\ddagger}$Department of Mathematics, Faculty of Science\\
King Mongkut's University of Technology Thonburi,\\
Thailand\\
$^{\ddagger}$Mathematics and Statistics with Applications (MaSA)\\
\small \tt Email: pawaton.kae@kmutt.ac.th; odilefava@gmail.com}

\date{}
\maketitle


%

\begin{abstract}
A dominating set of a graph $G=(V,E)$ is a vertex set $D$ such that every vertex in $V(G) \setminus D$ is adjacent to a vertex in $D$. The cardinality of a smallest dominating set of $D$ is called the domination number of $G$ and is denoted by $\gamma(G)$.  A vertex set $D$ is a $k$-isolating set of $G$ if $G - N_{G}[D]$ contains no $k$-cliques. The minimum cardinality of a $k$-isolating set of $G$ is called the $k$-isolation number of $G$ and is denoted by $\iota_{k}(G)$. Clearly, $\gamma(G) = \iota_{1}(G)$. A vertex set $I$ is  irredundant if, for every non-isolated vertex $v$ of $G[I]$, there exists a vertex $u$ in $V \setminus I$ such that $N_{G}(u) \cap I = \{v\}$. An irredundant set $I$ is maximal if the set $I \cup \{u\}$ is no longer irredundant for any $u \in V(G) \setminus I$. The minimum cardinality of a maximal irredundant set is called the irredundance number of $G$ and is denoted by $ir(G)$. Allan and Laskar \cite{AL1978} and Bollob\'{a}s and Cockayne \cite{BoCo1979} independently proved that $\gamma(G) < 2ir(G)$, which can be written $\iota _1(G) < 2ir(G)$, for any graph $G$.
In this paper, for a graph $G$ with maximum degree $\Delta$, we establish sharp upper bounds on $\iota_{k}(G)$ in terms of $ir(G)$ for $\Delta - 2 \leq k \leq \Delta + 1$.
\end{abstract}


\section{Notation}

Let $G$ be a simple graph with vertex set $V(G)$ and edge set $E(G)$ of order~$n(G) = |V(G)|$ and size $m(G) = |E(G)|$. We denote the \emph{degree} of $v$ in $G$ by $deg_G(v)$ and denote the maximum degree of $G$ by $\Delta(G)$. A \emph{neighbor} of a vertex $v$ in $G$ is a vertex $u$ which is adjacent to $v$. The \emph{open neighbor set} $N_G(v)$ of a vertex $v$ in $G$ is the set of neighbors of $v$. That is $N_G(v) = \{u \in V(G) \, | \, uv \in E(G)\}$. The \emph{closed neighbor set} of $v$ is $N_G[v] = N_G(v) \cup \{v\}$. For a subset $S \subseteq V(G)$, we use $N_{S}(v)$ to denote $N_{G}(v) \cap S$ and $deg_{S}(v) = |N_{G}(v) \cap S|$. If $v\in S$ and $deg_{S}(v) = 0$, then $v$ is an \emph{isolated vertex} in $S$. The \emph{neighbor set} of a vertex subset $S$ of $G$ is the set $N_G(S) = \cup_{v \in S} N_G(v)$. The \emph{closed neighbor set} of $S$ in $G$ is the set $N_G[S] = N_G(S) \cup S$. For a vertex $x \in S$, the \emph{private neighbor set} of $x$ with respect to $S$ is $N_{G}[x] \setminus N_{G}[S \setminus \{x\}] $ and is denoted by $PN(x, S)$. In particular, a vertex $y \in V(G) \setminus S$ is a \emph{private neighbor} of $x$ with respect to $S$ if $N_{S}(y) = \{x\}$. It is worth noting that if $x$ an isolated vertex in $S$, then $x \in PN(x, S)$. The subgraph of $G$ induced by $S$ is denoted by $G[S]$. The subgraph obtained from $G$ by deleting all vertices in $S$ and all edges incident with vertices in $S$ is denoted by $G - S$. The \emph{distance} between two vertices $u$ and $v$ in a connected graph $G$ is the length of a shortest $(u,v)$-path in $G$ and is denoted by $d_G(u,v)$. A \emph{complete graph} on $n$ vertices is denoted by $K_n$. For a graph $G$, if $H$ is a subgraph of $G$ isomorphic to $K_{k}$, then $H$ is called a $k$\emph{-clique}. A $k$-clique $H$ is \emph{adjacent to} a vertex $v$ (vice versa) if $v \in V(H)$ or $v$ is adjacent to a vertex of $H$. A \emph{block} $B$ of a graph $G$ is a maximal subgraph of $G$ such that $B$ itself does not contain a cut vertex. A graph $G$ is a \emph{block graph} if every block of $G$ is a complete graph while $G$ is a \emph{cactus} if each edge of $G$ belongs to at most one cycle. The \emph{cyclomatic number} of $G$ is $\mu(G) = m(G) - n(G) + k(G)$ where $k(G)$ is the number of components of $G$.
\vskip 5 pt

A subset $I \subseteq V(G)$ is called \emph{irredundant} if $PN(x, I) \neq \emptyset$ for all $x \in I$. Thus, for a vertex $x \in I$, there exists a vertex $y \in V(G) \setminus I$ such that $N_{I}(y) = \{x\}$ or $x$ is an isolated vertex in $I$. An irredundant set $I$ is said to be \emph{maximal} if $I$ is not a proper subset of any irredundant set of $G$. The cardinality of a smallest maximal irredundant set of $G$ is called the \emph{irredundance number} of $G$ and is denoted by $ir(G)$.
\vskip 5 pt

For vertex subsets $S$ and $R$ of $G$, we say that $S$ \emph{dominates} $R$ if every vertex in $R \setminus S$ is adjacent to a vertex in $S$. We write $S \succ R$ if $S$ dominates $R$. If $S = \{s\}$, we write $s \succ R$. Moreover, $S$ is a \emph{dominating set} of $G$ if $S \succ V(G)$. The cardinality of a smallest dominating set of $G$ is called the \emph{domination number} of $G$ and is denoted by $\gamma(G)$.

 Let $\cF$ be a family of graphs. A vertex subset $S \subseteq V(G)$ is said to be
$\cF$\emph{-isolating} if $G -  N_{G}[S]$ does not contain $H$ as an induced subgraph for all $H \in \cF$.
The $\cF$-\emph{isolation number} of $G$, denoted by $\iota_{\cF}(G)$, is the minimum cardinality of an $\cF$-isolating set of $G$.
This notion was introduced by Caro and Hansberg \cite{CaHa17}. If $\cF = \{K_{k}\}$, we say that a $\{K_k\}$-isolating set is $K_k$-isolating and we denote $\iota_{\{K_{k}\}}(G)$ by $\iota_{k}(G)$. A smallest $K_{k}$-isolating set is called a $\iota_{k}$-\emph{set}. Clearly, the $K_1$-isolating sets are the dominating sets and $\iota _1(G)=\gamma(G).$ Moreover every $K_k$-isolating set is $K_{k+1}$-isolating and $\iota_{k+1}(G)\le \iota_k(G)$.

Every minimal dominating set is irredundant and thus $ir(G)\le \gamma (G)$ for every graph $G$.
An upper bound on $\gamma(G)$ in terms of $ir(G)$ was independently established by Allan and Laskar \cite{AL1978} and by Bollob\'{a}s and Cockayne \cite{BoCo1979}.

\begin{thm}\label{thm CoBo}\cite{AL1978, BoCo1979}
Let $G$ be a graph. Then
$\gamma(G) < 2ir(G)$.
\end{thm}

Several authors lowered this bound for particular classes of graphs.
In 1991,  Damaschke \cite{D} reduced the upper bound of Theorem \ref{thm CoBo} for trees.

\begin{thm}\label{Dam}\cite{D}
Let $T$ be a tree. Then
$\gamma(T) < 3ir(T)/2$.
\end{thm}

\noindent The upper bound of Theorem \ref{Dam} was proved to be true by Volkman \cite{V}, in 1998, for some classes of tree-like graphs as block graphs, cacti and the graphs whose cyclomatic number is at most two. Volkman further conjectured that:

\begin{conj}\cite{V}
If $G$ is a cactus, then $5\gamma(G) < 8 ir(G)$.
\end{conj}

\noindent Volkman's conjecture was proved  by Zverovich \cite{Z}. Favaron, Kabanov and Puech \cite{FKP} showed that $\gamma(G) \leq 3 ir(G)/2$ if $G$ is claw-free and that the bound is sharp in different classes of claw-free graphs.

\indent Theorem 1 can be written $\iota_1(G) \le 2ir(G)-1$. Our aim  is to generalize this result and to find a sharp upper bound on $\iota _k(G)$ in terms of $ir(G)$ for some larger values of $k$.

\section{Main results}\label{sec2}

If a graph $G$ with maximum degree $\Delta(G) = \Delta$ contains a $k$-clique then
$k \leq \Delta + 1$ since every vertex in a $k$-clique is adjacent to at least $k - 1 \leq \Delta$ vertices.
When $\Delta \leq 1$, the graph $G$ is the union of $p_{1}$ paths of length one and $p_0$ isolated vertices
yielding that $p_{1} = \iota_2(G) \leq \iota_1(G)= p_0 + p_1=ir(G)$. We assume throughout that $\Delta \geq 2$ and we consider the four cases $\Delta - 2 \leq k \leq \Delta + 1$. In this section, we give the bounds on $\iota _k$ in terms of $ir$ without proofs and show that the bounds are attained. The proofs are given in the
following sections.

\begin{thm}\label{thm 1}
Let $G$ be a graph with maximum degree $\Delta \ge 2$. If $\Delta \leq k \leq \Delta + 1$, then $\iota_{k}(G) \leq ir(G)$.
\end{thm}
\vskip 5 pt

\indent To show the sharpness of the bound when $k = \Delta + 1$, we let $G^{1}_{t}$ be the disjoint union of $t$ copies $K^{1}_{k}, K^{2}_{k}, ..., K^{t}_{k}$ of a clique of $k = \Delta + 1$ vertices. Let $x_{i} \in V(K^{i}_{k})$. We see that $deg_{G^{1}_{t}}(x_{i}) = \Delta$. Clearly, $\{x_{1}, x_{2}, ..., x_{t}\}$ is an $\iota_{k}$-set and also a smallest maximal irredundant set. Therefore, $ir(G^{1}_{t}) = t = \iota_{k}(G^{1}_{t})$.

When $k = \Delta$, we let $K^{1}_{k}, K^{2}_{k}, ..., K^{2t}_{k}$ be $2t$ copies of a clique of $k = \Delta$ vertices. Let $x_{i} \in V(K^{i}_{k})$ and let $P^{i}_{2} = x^{i}_{1}x^{i}_{2}$ be $t$ copies of paths of $2$ vertices. The graph $G^{2}_{2t}$ is obtained by joining $x^{i}_{1}$ to $x_{2i - 1}$ and $x^{i}_{2}$ to $x_{2i}$. So, $deg_{G^{2}_{2t}}(x_{i}) = \Delta$ for all $1 \leq i \leq 2t$. Clearly, $\{x_{1}, x_{2}, ..., x_{2t}\}$ is an $\iota_{k}$-set, implying that $\iota_{k}(G^{2}_{2t}) = 2t$. Now, we let $I$ be a smallest maximal irredundant set of $G^{2}_{2t}$. If $I \cap (\{x^{i}_{1}\} \cup V(K^{2i-1}_{k})) = \emptyset$ for some $1 \leq i \leq t$, then $x'_{2i - 1}$ has a private neighbor with respect to $I \cup \{x'_{2i - 1}\}$ where $x'_{2i - 1} \in V(K^{2i-1}_{k}) \setminus \{x_{2i - 1}\}$. Thus, $I \cup \{x'_{2i - 1}\}$ is an irredundant set of $G$ containing $I$ contradicting the maximality of $I$. Hence, $I \cap (\{x^{i}_{1}\} \cup V(K^{2i-1}_{k})) \neq \emptyset$ for all $1 \leq i \leq t$. Similarly, $I \cap (\{x^{i}_{2}\} \cup V(K^{2i}_{k})) \neq \emptyset$. This yields $ir(G^{2}_{2t}) = |I| \geq 2t$. On the other hand, we see that $\{x_{1}, x_{2}, ..., x_{2t}\}$ is a maximal irredundant set of $G^{2}_{2t}$. By the minimality of $ir(G^{2}_{2t})$, we have $ir(G^{2}_{2t}) \leq 2t$. These implies that $ir(G^{2}_{2t}) = 2t = \iota_{k}(G^{2}_{2t})$.
\vskip 3 pt

\begin{thm}\label{thm 2}
Let $G$ be a graph with maximum degree $\Delta \ge 2$. If $k = \Delta - 1$, then $\iota_{k}(G) \leq \frac{(3\Delta - 4)ir(G)}{2\Delta - 2}$.
\end{thm}

\indent When $\Delta = k +1 \geq 7$, we show the sharpness of the bound of Theorem \ref{thm 2} by constructing connected graphs $G(k, l)$ that satisfy $\iota_{k}(G(k, l)) = \Big\lfloor\frac{(3\Delta - 4)ir(G(k, l))}{2\Delta - 2}\Big\rfloor$.  We first construct the graph $L(k)$ as follows. Let $x_{1}x_{2}x_{3}$ be a path of length two, $y_{1}, y_{2}, y_{3}$ be three distinct vertices and $K^{0}_{k}, K^{1}_{k}, K^{2}_{k}$ and $K^{3}_{k}$ be four disjoint $k$-cliques. The graph $L(k)$ is constructed as follows:
\begin{itemize}
  \item for $1 \leq i \leq 3$, join $y_{i}$ to all vertices in $\{x_{i}\} \cup V(K^{i}_{k})$,
  \item join $x_{2}$ to $\Delta - 3$ vertices in $K^{0}_{k}$ and, for $i \in \{1, 3\}$, join $x_{i}$ to $t_{i}$ vertices in $K^{0}_{k}$ in such a way that (\emph{i}) $t_{1} + t_{3} = \Delta + 1$ and $1 \leq t_{1}, t_{3} \leq \Delta - 3$ and (\emph{ii}) every vertex in $K^{0}_{k}$ is adjacent to exactly $2$ vertices in $\{x_{1}, x_{2}, x_{3}\}$.
\end{itemize}

\begin{figure}[H]
\centering
\definecolor{ududff}{rgb}{0.30196078431372547,0.30196078431372547,1}
\resizebox{1\textwidth}{!}{%
\begin{tikzpicture}[line cap=round,line join=round,>=triangle 45,x=1cm,y=1cm]

\draw [line width=0.5pt] (1,0)-- (3,0);
\draw [line width=0.5pt] (6,0)-- (8,0);
\draw [line width=0.5pt] (11,0)-- (13,0);
\draw [line width=0.5pt] (1,0)-- (1,-2);
\draw [line width=0.5pt] (1,0)-- (2,-2);
\draw [line width=0.5pt] (1,0)-- (1.48,-2);
\draw [line width=0.5pt] (1,0)-- (2.5,-2);
\draw [line width=0.5pt] (2,0)-- (3.5,-2);
\draw [line width=0.5pt] (2,0)-- (1,-2);
\draw [line width=0.5pt] (2,0)-- (3,-2);
\draw [line width=0.5pt] (2,0)-- (1.48,-2);
\draw [line width=0.5pt] (3,0)-- (2,-2);
\draw [line width=0.5pt] (3,0)-- (2.5,-2);
\draw [line width=0.5pt] (3,0)-- (3,-2);
\draw [line width=0.5pt] (3,0)-- (3.5,-2);
\draw [line width=0.5pt] (6,0)-- (5.755122099559129,-1.9981115760808803);
\draw [line width=0.5pt] (6,0)-- (6.258093295845571,-1.9981115760808803);
\draw [line width=0.5pt] (6,0)-- (6.770935622506731,-1.9976270136239493);
\draw [line width=0.5pt] (6,0)-- (7.2586478978749955,-1.9976270136239493);
\draw [line width=0.5pt] (7,0)-- (7.754355456446017,-1.9976270136239493);
\draw [line width=0.5pt] (7,0)-- (8.250063015017039,-1.9976270136239493);
\draw [line width=0.5pt] (7,0)-- (5.755122099559129,-1.9981115760808803);
\draw [line width=0.5pt] (7,0)-- (6.258093295845571,-1.9981115760808803);
\draw [line width=0.5pt] (8,0)-- (6.770935622506731,-1.9976270136239493);
\draw [line width=0.5pt] (8,0)-- (7.2586478978749955,-1.9976270136239493);
\draw [line width=0.5pt] (8,0)-- (7.754355456446017,-1.9976270136239493);
\draw [line width=0.5pt] (8,0)-- (8.250063015017039,-1.9976270136239493);
\draw [line width=0.5pt] (11,0)-- (10.5,-2);
\draw [line width=0.5pt] (11,0)-- (11,-2);
\draw [line width=0.5pt] (11,0)-- (11.5,-2);
\draw [line width=0.5pt] (11,0)-- (12,-2);
\draw [line width=0.5pt] (12,0)-- (12.5,-2);
\draw [line width=0.5pt] (12,0)-- (13,-2);
\draw [line width=0.5pt] (12,0)-- (10.5,-2);
\draw [line width=0.5pt] (12,0)-- (11,-2);
\draw [line width=0.5pt] (13,0)-- (11.5,-2);
\draw [line width=0.5pt] (13,0)-- (12,-2);
\draw [line width=0.5pt] (13,0)-- (12.5,-2);
\draw [line width=0.5pt] (13,0)-- (13,-2);
\draw [rotate around={-0.04651406773445012:(2.3016722224108146,-1.9986212570749933)},line width=0.6pt] (2.3016722224108146,-1.9986212570749933) ellipse (1.793396780179797cm and 0.5761535126128963cm);
\draw [rotate around={-0.1796100057470778:(7.001773068267073,-2.00258790147116)},line width=0.6pt] (7.001773068267073,-2.00258790147116) ellipse (1.798703267426045cm and 0.5724255105264566cm);
\draw [rotate around={-0.12960621793887145:(11.753554719067996,-2.003966644396167)},line width=0.6pt] (11.753554719067996,-2.003966644396167) ellipse (1.8400925961586179cm and 0.5576476265568429cm);
\draw [line width=0.5pt] (1,3) circle (0.4008235990407263cm);
\draw [line width=0.5pt] (2,4) circle (0.4024322523463133cm);
\draw [line width=0.5pt] (3,3) circle (0.393810397047737cm);
\draw [line width=0.5pt] (6,3) circle (0.38339920023310353cm);
\draw [line width=0.5pt] (7,4) circle (0.4023850267161873cm);
\draw [line width=0.5pt] (8,3) circle (0.39773414025041776cm);
\draw [line width=0.5pt] (11,3) circle (0.38325629187001375cm);
\draw [line width=0.5pt] (12,4) circle (0.3898946807091898cm);
\draw [line width=0.5pt] (13,3) circle (0.38703549457956476cm);
\draw [line width=0.5pt] (1,0)-- (1,1);
\draw [line width=0.5pt] (2,0)-- (2,2);
\draw [line width=0.5pt] (3,0)-- (3,1);
\draw [line width=0.5pt] (6,0)-- (6,1);
\draw [line width=0.5pt] (7,0)-- (7,2);
\draw [line width=0.5pt] (8,0)-- (8,1);
\draw [line width=0.5pt] (11,0)-- (11,1);
\draw [line width=0.5pt] (12,0)-- (12,2);
\draw [line width=0.5pt] (13,0)-- (13,1);
\draw [line width=0.5pt] (3,0)-- (6,0);
\draw [line width=0.5pt] (8,0)-- (11,0);
\draw [line width=0.5pt] (1,1)-- (0.6437084157120798,2.8163811978212343);
\draw [line width=0.5pt] (1,1)-- (1.3697855743254315,2.8453449432809528);
\draw [line width=0.5pt] (2,2)-- (1.6359382530372755,3.828504337877088);
\draw [line width=0.5pt] (2,2)-- (2.3510318224443756,3.803204732374208);
\draw [line width=0.5pt] (3,1)-- (2.655950451579013,2.8083844028942027);
\draw [line width=0.5pt] (3,1)-- (3.3461608426571927,2.812225401521509);
\draw [line width=0.5pt] (6,1)-- (5.646528902159162,2.85149703797373);
\draw [line width=0.5pt] (6,1)-- (6.355619601335921,2.8567182988478157);
\draw [line width=0.5pt] (7,2)-- (6.644959427955781,3.810632891113845);
\draw [line width=0.5pt] (7,2)-- (7.363360449732776,3.827133307732861);
\draw [line width=0.5pt] (8,1)-- (7.63072461280423,2.8522565238930895);
\draw [line width=0.5pt] (8,1)-- (8.362603450067276,2.836564433735961);
\draw [line width=0.5pt] (11,1)-- (10.647309209015532,2.8500180303746387);
\draw [line width=0.5pt] (11,1)-- (11.353916814834646,2.8529344587088743);
\draw [line width=0.5pt] (12,2)-- (11.656207163875527,3.816086031322217);
\draw [line width=0.5pt] (12,2)-- (12.348874989408667,3.8259192606563737);
\draw [line width=0.5pt] (13,1)-- (12.644682061547035,2.8465541245976844);
\draw [line width=0.5pt] (13,1)-- (13.339800620225532,2.8147110025965065);
\begin{scriptsize}
\draw [fill=black] (1,0) circle (2pt);
\draw[color=black] (1.2,0.2) node {$x^{1}_{1}$};
\draw [fill=black] (2,0) circle (2pt);
\draw[color=black] (2.2,0.2) node {$x^{1}_{2}$};
\draw [fill=black] (11,0) circle (2pt);
\draw[color=black] (11.2,0.2) node {$x^{3}_{1}$};
\draw [fill=black] (13,0) circle (2pt);
\draw[color=black] (13.2,0.2) node {$x^{3}_{3}$};
\draw [fill=black] (12,0) circle (2pt);
\draw[color=black] (12.2,0.2) node {$x^{3}_{2}$};
\draw [fill=black] (8,0) circle (2pt);
\draw[color=black] (8.2,0.2) node {$x^{2}_{3}$};
\draw [fill=black] (7,0) circle (2pt);
\draw[color=black] (7.2,0.2) node {$x^{2}_{2}$};
\draw [fill=black] (6,0) circle (2pt);
\draw[color=black] (6.2,0.2) node {$x^{2}_{1}$};
\draw [fill=black] (1,-2) circle (2pt);
\draw [fill=black] (2,-2) circle (2pt);
\draw [fill=black] (3,-2) circle (2pt);
\draw [fill=black] (13,-2) circle (2pt);
\draw [fill=black] (12,-2) circle (2pt);
\draw [fill=black] (11,-2) circle (2pt);
\draw [fill=black] (5.755122099559129,-1.9981115760808803) circle (2pt);
\draw [fill=black] (6.770935622506731,-1.9976270136239493) circle (2pt);
\draw [fill=black] (7.754355456446017,-1.9976270136239493) circle (2pt);
\draw [fill=black] (3,0) circle (2pt);
\draw[color=black] (3.2,0.2) node {$x^{1}_{3}$};
\draw [fill=black] (1.48,-2) circle (2pt);
\draw [fill=black] (2.5,-2) circle (2pt);
\draw [fill=black] (3.5,-2) circle (2pt);
\draw [fill=black] (6.258093295845571,-1.9981115760808803) circle (2pt);
\draw [fill=black] (7.2586478978749955,-1.9976270136239493) circle (2pt);
\draw [fill=black] (8.250063015017039,-1.9976270136239493) circle (2pt);
\draw [fill=black] (10.5,-2) circle (2pt);
\draw [fill=black] (11.5,-2) circle (2pt);
\draw [fill=black] (12.5,-2) circle (2pt);
\draw [fill=black] (1,1) circle (2pt);
\draw[color=black] (1.3,1.2) node {$y^{1}_{1}$};
\draw [fill=black] (2,2) circle (2pt);
\draw[color=black] (2.3,2.2) node {$y^{1}_{2}$};
\draw [fill=black] (3,1) circle (2pt);
\draw[color=black] (3.3,1.2) node {$y^{1}_{3}$};
\draw [fill=black] (6,1) circle (2pt);
\draw[color=black] (6.3,1.2) node {$y^{2}_{1}$};
\draw [fill=black] (7,2) circle (2pt);
\draw[color=black] (7.3,2.2) node {$y^{2}_{2}$};
\draw [fill=black] (8,1) circle (2pt);
\draw[color=black] (8.3,1.2) node {$y^{2}_{3}$};
\draw [fill=black] (11,1) circle (2pt);
\draw[color=black] (11.3,1.2) node {$y^{3}_{1}$};
\draw [fill=black] (12,2) circle (2pt);
\draw[color=black] (12.3,2.2) node {$y^{3}_{2}$};
\draw [fill=black] (13,1) circle (2pt);
\draw[color=black] (13.3,1.2) node {$y^{3}_{3}$};

\end{scriptsize}
\end{tikzpicture}
 }%
\vskip -0.25 cm
\caption{The graph $G(6, 3)$.}
\label{g1}
\end{figure}
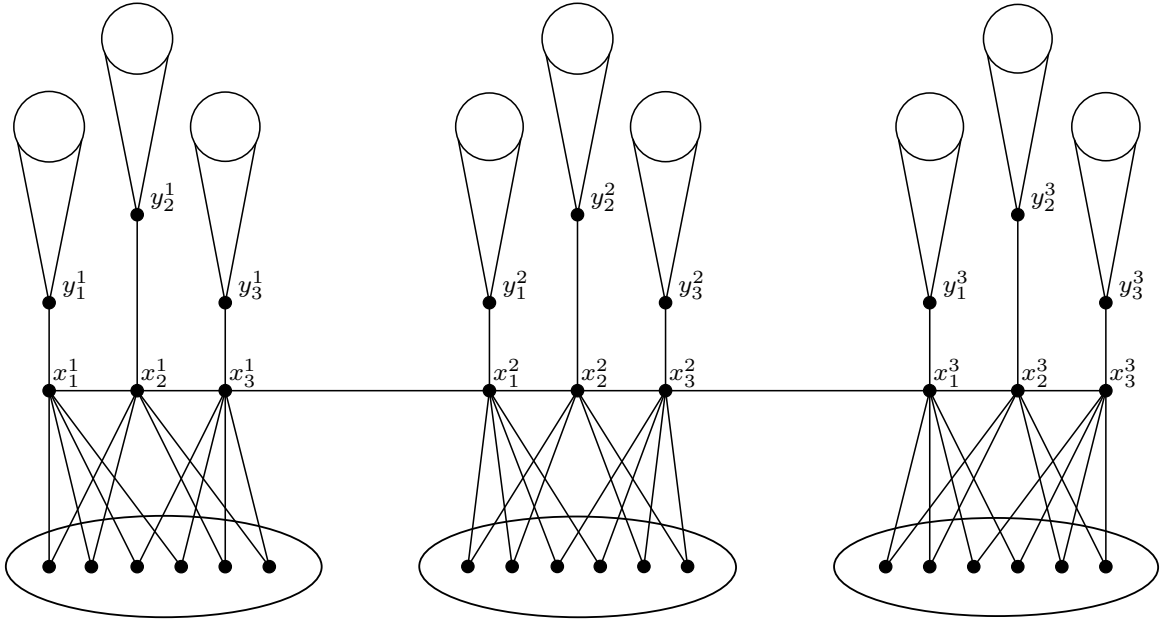

\noindent The above construction is possible because there are $2k$ edges from $K^{0}_{k}$ to $\{x_{1}, x_{2}, x_{3}\}$ and there are $\Delta - 3 + \Delta + 1 = 2\Delta - 2 = 2k$ edges from $\{x_{1}, x_{2}, x_{3}\}$ to $K^{0}_{k}$. It can be observed that the vertex $y_{i}$ has degree $k + 1 = \Delta$ and the vertex $x_{2}$ is adjacent to $\Delta - 3$ vertices in $K^{0}_{k}$. Further, $deg_{L(k)}(x_{1}) = deg_{L(k)}(x_{3}) \leq \Delta - 1$ while $deg_{L(k)}(x_{2}) = \Delta$. Now, we let $l$ be a positive integer such that
\begin{center}
$4l = \lfloor\frac{(3\Delta - 4)3l}{2\Delta - 2}\rfloor$.
\end{center}
Note that $\ell$ is at most 3 for $\Delta =7$, at most 2 until $\Delta =10$ and $\ell=1$ for bigger $\Delta$.
\noindent The graph $G(k, l)$ that satisfies $\iota_{k}(G(k, l)) = \Big\lfloor\frac{(3\Delta - 4)ir(G(k, l))}{2\Delta - 2}\Big\rfloor$ is obtained from $l$ copies $L_{1}, ..., L_{l}$ of the graph $L(k)$ by joining $x_{3}$ of $L_{i}$ to $x_{1}$ of $L_{i + 1}$ for all $1 \leq i \leq l - 1$. For the sake of convenient, we may relabel $x_{1}, x_{2}, x_{3}, y_{1}, y_{2}, y_{3}, K^{0}_{k}, K^{1}_{k}, K^{2}_{k}$ and $K^{3}_{k}$ in $L(k)$ to be $x^{i}_{1}, x^{i}_{2}, x^{i}_{3}, y^{i}_{1}, y^{i}_{2}, y^{i}_{3}, K^{i, 0}_{k}, K^{i, 1}_{k}, K^{i, 2}_{k}$ and $K^{i, 3}_{k}$ respectively in $L_{i}$. We let $z_{i}$ be a vertex in $K^{i, 0}_{k}$. Clearly, $\cup^{l}_{i = 1}\{y^{i}_{1}, y^{i}_{2}, y^{i}_{3}, z_{i}\}$ is  a $K_k-$isolating set of $G(k, l)$. By the minimality of $\iota_{k}(G(k, l))$, we have $\iota_{k}(G(k, l)) \leq 4l$. Let $S$ be an $\iota_{k}$-set of $G(k, l)$. To be adjacent to $K^{i, 0}_{k}, K^{i, 1}_{k}, K^{i, 2}_{k}$ and $K^{i, 3}_{k}$ we have that $S \cap (y^{i}_{j} \cup V(K^{i, j}_{k})) \neq \emptyset$ and $S \cap (V(K^{i, 0}_{k}) \cup \{x^{i}_{1}, x^{i}_{2}, x^{i}_{3}\}) \neq \emptyset$ for all $1 \leq i \leq l$ and $1 \leq j \leq 3$. Therefore, $\iota_{k}(G(k, l)) \geq 4l$ implying that $\iota_{k}(G(k, l)) = 4l$. Next, we let $I$ be a smallest maximal irredundant set of $G(k, l)$. By maximality of $I$, we have that $I \cap (\{x^{i}_{j}, y^{i}_{j}\} \cup V(K^{i, j}_{k})) \neq \emptyset$. This yields $ir(G(k, l)) \geq 3l$. Clearly, $\cup^{l}_{i = 1}\{x^{i}_{1}, x^{i}_{2}, x^{i}_{3}\}$ is a maximal irredundant set of $G(k, l)$. Hence,
\begin{center}
$\iota_{k}(G(k, l)) = 4l = \lfloor\frac{(3\Delta - 4)3l}{2\Delta - 2}\rfloor = \lfloor\frac{(3\Delta - 4)ir(G(k, l))}{2\Delta - 2}\rfloor$.
\end{center}
\noindent An example of the graph $G(k, l)$ when
$k = 6, \Delta = k + 1 = 7$ and $l = 3$ is illustrated by Figure \ref{g1}. In the figure, we let a circle or an oval to denote a clique of order $6$. We see that $\iota_{6}(G(6, 3)) = 12$ while $ir(G(6, 3)) = 9$.
\vskip 18 pt

\vskip 18 pt

\begin{thm}\label{thm 3}
Let $G$ be a graph with maximum degree $\Delta \ge 3$. If $k = \Delta - 2$, then $\iota_{k}(G) \leq \frac{3ir(G)}{2}$.
\end{thm}

\indent For  connected graphs satisfying the upper bound of Theorem \ref{thm 3}, we let $x_{1}x_{2}x_{3}x_{4}$ be a path of length $3$, $K^{0}_{k}, K^{1}_{k}$ and $K^{2}_{k}$ be $3$ $k$-cliques. The graph $D(k)$ is constructed as follows:
\begin{itemize}
  \item join $x_{1}$ to all vertices in $V(K^{1}_{k})$,
  \item join $x_{4}$ to all vertices in $V(K^{2}_{k})$ and
  \item join both $x_{2}$ and $x_{3}$ to all vertices in $V(K^{0}_{k})$.
\end{itemize}
\noindent Then, for $t\ge 2$, the graph $D(k, t)$ is obtained by the disjoint union of $t$ copies $D_{1}, ..., D_{t}$ of $D(k)$. We relabel $x_{i}$ and $K^{j}_{k}$ of $D(k)$ to be $x^{\ell}_{i}$ and $K^{\ell, j}_{k}$ in $D_{\ell}$ respectively for all $1 \leq i \leq 4,~ 0 \leq j \leq 2$ and $1 \leq \ell \leq t$.  Then, add the edges $x_4^{\ell}x_1^{\ell +1}$ for $1 \leq \ell \leq t - 1$.
It can be easily checked  that
$\cup^{t}_{\ell = 1}\{ x^{\ell}_{1}, x^{\ell}_{2}, x^{\ell}_{4}\}$
is a smallest $\iota_{k}$-set and $\cup^{t}_{\ell = 1}\{x^{\ell}_{2}, x^{\ell}_{3}\}$ is a smallest maximal irredundant set of $D(k, t)$. Thus $\iota_{k}(D(k, t)) = 3t$ and $ir(D(k, t)) = 2t$ implying that $\iota_{k}(D(k, t)) = \frac{3ir(D(k, t))}{2}$.

\vskip 5 pt

We saw that for $k=1$, $\iota_1(G)=\gamma(G)$. Hence, Theorem \ref{thm 2} confirms the fact that paths and cycles have equal irredundance and domination numbers and Theorem \ref{thm 3} admits the following Corollary.

\begin{cor} \label{cor 1}
Let $G$ be a graph with maximum degree $\Delta = 3$. Then  $\gamma(G)\le \frac{3ir(G)}{2}$ and the bound is sharp.
\end{cor}
Note that for cubic graphs, the slightly larger bound $\gamma (G)\le \frac{27}{16}ir(G)$ can be deduced from the properties $\gamma (G)\le \frac{3n}{8}$ \cite{Re} and $ir(G)\ge \frac{2n}{3\Delta}$ \cite{CoMy}.

The following graph $H$ is an example of subcubic graph for which the equality in Corollary \ref{cor 1} is attained. $H$ consists of $k$ paths $x_1^ix_2^i \cdots x_7^i$ with $1\le i\le k$ and the $2k$ edges $x_3^ix_5^i$ and  $x_6^ix_2^{i+1 \,mod (k)}$. It can be checked that $\bigcup_{i=1}^k\{x_2^i,x_4^i,x_6^i\}$ is a $\gamma(H)$-set and $\bigcup_{i=1}^k\{x_3^i,x_5^i\}$ is an $ir(H)$. Hence $\gamma(H)= \frac{3ir(H)}{2}$.

\indent  Let $c_{\Delta}(k)$ be the least upper bound of $\frac{\iota_{k}(G)}{ir(G)}$ for all graphs $G$ of maximum degree $\Delta$. From our results in Theorems \ref{thm 1}, \ref{thm 2} and \ref{thm 3}, we proved that the graphs $G$ of given maximum degree $\Delta$ satisfy $\iota_{k}(G) \leq c_{\Delta}(k)ir(G)$ for $\Delta -2\le k\le \Delta + 1$ with
 $c_{\Delta}(\Delta -2)=\frac{3}{2}$,  $c_{\Delta}(\Delta -1)=\frac{3\Delta -4}{2\Delta -2}$ and $c_{\Delta}(\Delta + 1)=c_{\Delta}(\Delta)=1$. For the smallest possible value of $k$ which is $1$, it follows by Theorem \ref{thm CoBo} that $c_{\Delta}(1) < 2$. For an arbitrary large $\Delta$, we conjecture that:
\vskip 2 pt

\begin{conj}
Let $G$ be a graph with maximum degree $\Delta$. Then $\iota_{k}(G) \leq c_{\Delta}(k)ir(G)$ where $c_{\Delta}(1), c_{\Delta}(2), ..., c_{\Delta}(\Delta)$ is a strictly decreasing sequence such that $1 \leq c_{\Delta}(k) < 2$ for  $1 \leq k \leq \Delta$.
\end{conj}
\vskip 3 pt

\noindent Now, we fix $k=1$ and let $c_{\Delta}(1)=c_{\Delta}$. By Corollary \ref{cor 1}, $c_{3} = \frac{3}{2}$.
We think that $c_{\Delta}$ converges to $2$ when $\Delta \rightarrow \infty$ and we we give a second conjecture generalizing Corollary \ref{cor 1}.

\begin{conj}
Let $G$ be a graph with maximum degree $\Delta \ge 2$. Then $\gamma(G) \leq c_{\Delta}ir(G)$ where $c_{\Delta}$ is a strictly increasing sequence such that $c_2=1$ and  $c_{\Delta} \rightarrow 2$ when $\Delta \rightarrow \infty$.
\end{conj}

The following example shows that if the previous conjecture is true, then $c_4\ge \frac{5}{3}=\frac{2\Delta -3}{\Delta -1}$.
The graph $G$ is constructed from 3 disjoint paths $P_6$, $c_{1i}b_{1i}a_{1i}a_{2i}b_{2i}c_{2i}$, $1\le i\le 3$ and six new vertices $u,v,w,x,y,z$ by adding the 12 edges $ua_{11}$, $ua_{12}$, $va_{21}$, $va_{22}$, $wa_{11}$, $wa_{13}$, $xa_{12}$, $xa_{13}$, $ya_{21}$, $ya_{23}$, $za_{22}$, $za_{23}$. Then
$\bigcup_{1\le i\le 3}\{a_{1i}, a_{2i}\}$ is a smallest maximal irredundant set and
$\bigcup_{1\le i\le 3}\{b_{1i}, b_{2i}\}\cup \{a_{21},a_{13},u,z\}$ is a minimum dominating set.
Thus $ir(G)=6$, $\gamma(G) =10$ and
$\frac{\gamma(G)}{ir(G)}=\frac{5}{3}=\frac{2\Delta -3}{\Delta -1}$.

\section{Preliminary and preparation}
\label{S:known}

We need the following properties of a maximal irredundant set which were established by Bollob$\acute{a}$s and Cockayne~\cite{BoCo1979}.

\begin{thm}\label{thm BoCo2}~\cite{BoCo1979}
Let $I$ be a maximal irredundant set of a graph $G$. If there exists a vertex $u$ which is not dominated by any vertex in $I$, then there exists a vertex $x \in I$ such that
\begin{enumerate}
  \item $PN(x, I) \subseteq N(u)$ and
  \item if there exists a pair of non-adjacent vertices
$x_1$ and $x_2$ in $PN(x, I)$, then, for each $1 \leq i \leq 2$, there exists $y_{i} \in I \setminus \{x\}$ such that $x_{i} \succ PN(y_{i}, I)$.
\end{enumerate}
\end{thm}

\noindent
\textbf{Observations:} \\
\noindent 1. In part (a) of Theorem \ref{thm BoCo2}, the vertices $x$ of $I$ such that
$PN(x, I) \subseteq N(u)$ are not isolated in $I$ since in this case, $x$ would belong to $PN(x,I)$ yielding that $PN(x,I)$ could not be dominated by $u$.

\noindent 2. Similarly in (b) the vertex $y_i$ is not isolated in $I$ for otherwise $y_i\in$ PN$(y_i,I)$ and the private neighbor $x_i$ of $x$ would be adjacent to another vertex $y_i$ of $X$, a contradiction.



We give now the notation that we will use  in the proofs of our main results. In the following, let $G$ be a graph. We let

\indent $I =$ a smallest maximal irredundant set of $G$,

\indent $U =$ the set of vertices which are not dominated by $I$,

\indent $P =$ the set of vertices $y \in PN(x, I) \setminus I$ for some $x \in I$ and

\indent $S = N(I) \setminus (P \cup I)$.

\noindent Therefore, every vertex in $S$ is adjacent to at least two vertices in $I$.  Clearly, the sets $I, U, P$ and $S$ form a partition of $V(G)$. By Theorem \ref{thm BoCo2}(a), for every vertex $u \in U$, there exists a non-isolated vertex $x \in I$ such that $PN(x, I) \subseteq N(u)$. Further, we partition the set $I$  as follows.

\indent $Z =$ the set of all isolated vertices of $I$,

\indent $A = I \setminus Z$, the set of non isolated vertices of $I$,

\indent $Q =$ the set of vertices $q$ in $A$ such that $PN(q, I) \nsubseteq N(u)$ for all $u \in U$,

\indent $B = A \setminus Q$, the set of non isolated vertices $q$ of $I$ whose $I$-private neighborhood is dominated by some $u\in U$,

\indent $N =$ a maximal set of vertices $x$ in $B$ such that $PN(x, I)$ is a clique or, for each $x' \in PN(x, I)$, $x' \succ PN(y, I)$ for some $y \in N$,

\noindent and

\indent $M = B \setminus N$.

\noindent  We note that if $U\neq \emptyset$, then $B \neq \emptyset$ by Theorem \ref{thm BoCo2}(a) and the partition of $B$ refers to
Theorem \ref{thm BoCo2}(b). From the definition, if $x \in M$ then $PN(x, I)$ is not a clique and thus $|PN(x, I)| \geq 2$, and  there exists a vertex $x' \in PN(x, I)$ such that
$x'$ does not dominate $PN(y, I)$ for any $y \in N$ but
$x' \succ PN(w, I)$ for some $w \in I \setminus N$ (by Theorem \ref{thm BoCo2}(b)).

Note also that for no vertex $z$ of $Z$, $PN(z,I)\subseteq N(u)$ since evry isolated vertex $z$ of $I$ belongs to $PN(z,I)$.

\begin{figure}[H]
\centering
\definecolor{ududff}{rgb}{0.30196078431372547,0.30196078431372547,1}
\resizebox{1.1\textwidth}{!}{%
\begin{tikzpicture}[line cap=round,line join=round,>=triangle 45,x=1cm,y=1cm]

\clip(-0.7877175596292759,-0.54661741023467) rectangle (13.489939214647773,9.547776869668256);
\draw [line width=0.7pt] (0,4)-- (0,3);
\draw [line width=0.7pt] (12,3)-- (0,3);
\draw [line width=0.7pt] (0,4)-- (12,4);
\draw [line width=0.7pt] (12,3)-- (12,4);
\draw [line width=0.7pt] (1,2)-- (1,1);
\draw [line width=0.7pt] (9,1)-- (1,1);
\draw [line width=0.7pt] (1,2)-- (9,2);
\draw [line width=0.7pt] (9,1)-- (9,2);
\draw [line width=0.7pt] (-0.02,4.58)-- (-0.02,6.38);
\draw [line width=0.7pt] (10.56,4.64)-- (10.58,6.38);
\draw [line width=0.7pt] (10.58,6.38)-- (-0.02,6.38);
\draw [line width=0.7pt] (-0.02,4.58)-- (10.56,4.64);
\draw [line width=0.7pt] (0,8)-- (0,7);
\draw [line width=0.7pt] (9.59454545454546,4)-- (9.59454545454546,3);
\draw [line width=0.7pt] (7,4)-- (7,3);
\draw [line width=0.7pt] (4.790054311816081,4)-- (4.79005431181608,3);

\draw [line width=0.7pt] (2.30471,3.40288)-- (3.39455,3.39818);
\draw [line width=0.7pt] (4.4,3.4)-- (6.2,3.4);
\draw [line width=0.7pt] (7.50143,3.40605)-- (8.31275,3.39956);
\draw [line width=0.7pt] (9.1046,3.39956)-- (8.31275,3.39956);
\draw [line width=0.7pt] (0.41273,3.39818)-- (1.5,1.5);
\draw [line width=0.7pt] (1.5,1.5)-- (3.39455,3.39818);
\draw [line width=0.7pt] (2.30471,3.40288)-- (3,1.5);
\draw [line width=0.7pt] (4.4,3.4)-- (5.5,1.5);
\draw [line width=0.7pt] (5.5,1.5)-- (9.90073,3.39665);
\draw [line width=0.7pt] (9.90073,3.39665)-- (9.89273,5.30381);
\draw [line width=0.7pt] (10.40297,3.40167)-- (10.3943,5.31077);
\draw [line width=0.7pt] (9.90073,3.39665)-- (9.5,5);
\draw [line width=0.7pt] (9.1046,3.39956)-- (9.10179,5.29839);
\draw [line width=0.7pt] (8.31275,3.39956)-- (8.3038,5.30539);
\draw [line width=0.7pt] (7.50143,3.40605)-- (7.5058,5.29839);
\draw [line width=0.7pt] (7.50143,3.40605)-- (7.8068,4.9064);
\draw [line width=0.7pt] (8.31275,3.39956)-- (8.61179,4.9204);
\draw [line width=0.7pt] (6.2,3.4)-- (6.18935,5.99099);
\draw [line width=0.7pt] (6.2,3.4)-- (6.57799,5.00171);
\draw [line width=0.7pt] (0.41273,3.39818)-- (0.20654,5.80256);
\draw [line width=0.7pt] (0.41273,3.39818)-- (0.59519,5.80256);
\draw [line width=0.7pt] (0.41273,3.39818)-- (0.40675,6.16765);
\draw [line width=0.7pt] (0.40675,6.16765)-- (0.20654,5.80256);
\draw [line width=0.7pt] (0.59519,5.80256)-- (0.40675,6.16765);
\draw [line width=0.7pt] (0.20654,5.80256)-- (0.59519,5.80256);
\draw [line width=0.7pt] (2.30471,3.40288)-- (2.29939,5.98948);
\draw [line width=0.7pt] (2.30471,3.40288)-- (2.5,5.5);
\draw [line width=0.7pt] (3.39455,3.39818)-- (3.40098,5.9955);
\draw [line width=0.7pt] (3.39455,3.39818)-- (3.20233,5.49587);
\draw [line width=0.7pt] (4.4,3.4)-- (4.38819,5.98948);
\draw [line width=0.7pt] (4.4,3.4)-- (4.19556,5.48986);
\draw [line width=0.7pt] (7,8)-- (7,7);
\draw [line width=0.7pt] (0,8)-- (7,8);
\draw [line width=0.7pt] (0,7)-- (7,7);
\draw [line width=0.7pt] (2.29939,5.98948)-- (3.40098,5.9955);
\draw [line width=0.7pt] (2.29939,5.98948)-- (3.20233,5.49587);
\draw [line width=0.7pt] (2.5,5.5)-- (3.40098,5.9955);
\draw [line width=0.7pt] (2.5,5.5)-- (3.20233,5.49587);
\draw [line width=0.7pt] (4.38819,5.98948)-- (3.40098,5.9955);
\draw [line width=0.7pt] (4.38819,5.98948)-- (3.20233,5.49587);
\draw [line width=0.7pt] (4.19556,5.48986)-- (3.40098,5.9955);
\draw [line width=0.7pt] (3.20233,5.49587)-- (4.19556,5.48986);
\draw [line width=0.7pt] (0.41273,3.39818)-- (2.30471,3.40288);
\draw [line width=0.7pt] (3,1.5)-- (0.41273,3.39818);
\draw [line width=0.7pt] (6.18935,5.99099)-- (9.10179,5.29839);
\draw [line width=0.7pt] (6.57799,5.00171)-- (7.5058,5.29839);
\draw [line width=0.7pt] (6.57799,5.00171)-- (7.8068,4.9064);
\draw [line width=0.7pt] (6.57799,5.00171)-- (4.19556,5.48986);
\draw [line width=0.7pt] (6.57799,5.00171)-- (4.38819,5.98948);
\draw [line width=0.7pt] (0.40176,7.49734)-- (0.40675,6.16765);
\draw [line width=0.7pt] (0.40176,7.49734)-- (0.20654,5.80256);
\draw [line width=0.7pt] (0.40176,7.49734)-- (0.59519,5.80256);
\draw [line width=0.7pt] (2.40044,7.49734)-- (2.29939,5.98948);
\draw [line width=0.7pt] (2.40044,7.49734)-- (2.5,5.5);
\draw [line width=0.7pt] (3.5,7.5)-- (4.38819,5.98948);
\draw [line width=0.7pt] (3.5,7.5)-- (3.40098,5.9955);
\draw [line width=0.7pt] (3.5,7.5)-- (3.20233,5.49587);
\draw [line width=0.7pt] (3.5,7.5)-- (4.19556,5.48986);
\draw [line width=0.7pt] (6.5,7.5)-- (6.18935,5.99099);
\draw [line width=0.7pt] (6.5,7.5)-- (6.57799,5.00171);
\begin{scriptsize}
\draw [fill=black] (0.41273,3.39818) circle (1.8pt);
\draw [fill=black] (2.30471,3.40288) circle (1.8pt);
\draw [fill=black] (3.39455,3.39818) circle (1.8pt);
\draw [fill=black] (4.4,3.4) circle (1.8pt);
\draw [fill=black] (6.2,3.4) circle (1.8pt);
\draw [fill=black] (7.50143,3.40605) circle (1.8pt);
\draw [fill=black] (8.31275,3.39956) circle (1.8pt);
\draw [fill=black] (9.1046,3.39956) circle (1.8pt);
\draw [fill=black] (9.90073,3.39665) circle (1.8pt);
\draw [fill=black] (10.40297,3.40167) circle (1.8pt);
\draw [fill=black] (11.00064,3.40167) circle (1.8pt);
\draw [fill=black] (11.5079,3.40167) circle (1.8pt);
\draw [fill=black] (1.5,1.5) circle (1.8pt);
\draw [fill=black] (3,1.5) circle (1.8pt);
\draw [fill=black] (5.5,1.5) circle (1.8pt);
\draw [fill=black] (9.89273,5.30381) circle (1.8pt);
\draw [fill=black] (10.3943,5.31077) circle (1.8pt);
\draw [fill=black] (9.5,5) circle (1.8pt);
\draw [fill=black] (9.10179,5.29839) circle (1.8pt);
\draw [fill=black] (8.3038,5.30539) circle (1.8pt);
\draw [fill=black] (7.5058,5.29839) circle (1.8pt);
\draw [fill=black] (7.8068,4.9064) circle (1.8pt);
\draw [fill=black] (8.61179,4.9204) circle (1.8pt);
\draw [fill=black] (6.18935,5.99099) circle (1.8pt);
\draw [fill=black] (6.57799,5.00171) circle (1.8pt);
\draw [fill=black] (0.20654,5.80256) circle (1.8pt);
\draw [fill=black] (0.59519,5.80256) circle (1.8pt);
\draw [fill=black] (0.40675,6.16765) circle (1.8pt);
\draw [fill=black] (2.29939,5.98948) circle (1.8pt);
\draw [fill=black] (3.40098,5.9955) circle (1.8pt);
\draw [fill=black] (3.20233,5.49587) circle (1.8pt);
\draw [fill=black] (4.38819,5.98948) circle (1.8pt);
\draw [fill=black] (4.19556,5.48986) circle (1.8pt);
\draw [fill=black] (0.40176,7.49734) circle (1.8pt);
\draw [fill=black] (2.5,5.5) circle (1.8pt);
\draw [fill=black] (2.40044,7.49734) circle (1.8pt);
\draw [fill=black] (3.5,7.5) circle (1.8pt);
\draw [fill=black] (6.5,7.5) circle (1.8pt);

\draw[color=black] (12.5,3.5) node {$I$};
\draw[color=black] (9.6,1.5) node {$S$};
\draw[color=black] (11,5.2) node {$P$};
\draw[color=black] (7.5,7.5) node {$U$};

\draw[color=black] (3,2.7) node {$N$};
\draw[color=black] (6,2.7) node {$M$};
\draw[color=black] (8,2.7) node {$Q$};
\draw[color=black] (11,2.7) node {$Z$};

\end{scriptsize}
\end{tikzpicture}
 }%
\vskip -1 cm
\caption{The partition $I \cup P \cup S \cup U$ of the graph and the partition $N \cup M \cup Q \cup Z$ of $I$.}
\label{ipsu}
\end{figure}
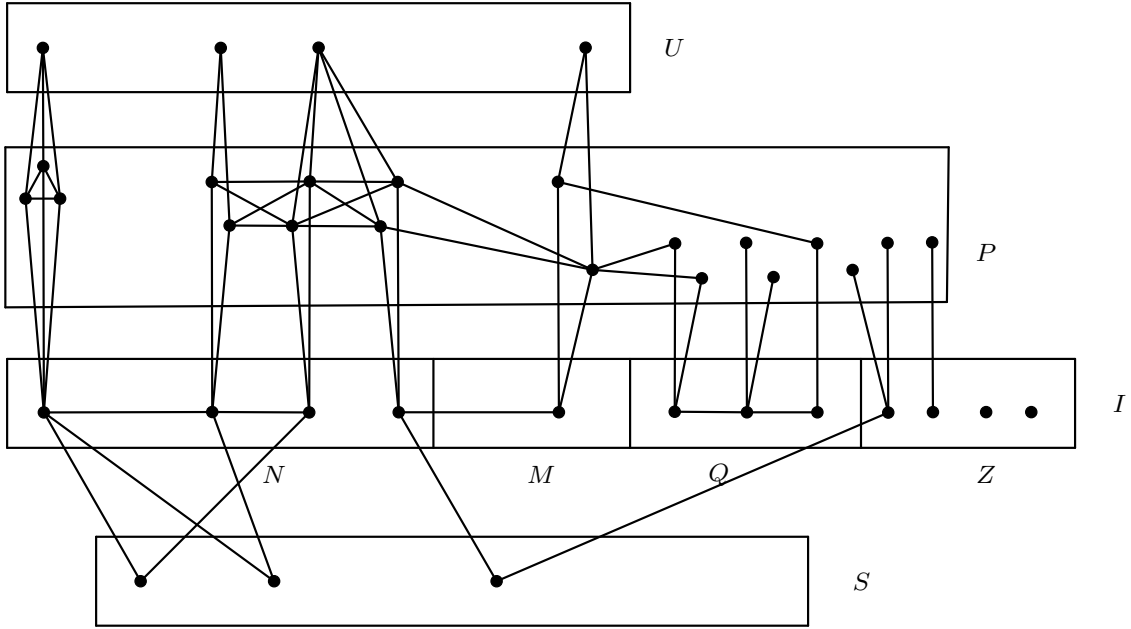

\section{Proofs}\label{S:main}




Throughout this section,
$I$ is a minimum maximal irredundant set as in Section 3.
Let $N = \{a_{1}, a_{2}, ..., a_{|N|}\}$ and choose $a'_{i} \in PN(a_{i}, I)$ for all $1 \leq i \leq |N|$.
The following lemma is a property of
$\hat{N} = \{a'_{1}, a'_{2}, ..., a'_{|N|}\}$.

\begin{lem}\label{l1}
$\hat{N} \succ \cup^{|N|}_{i = 1}PN(a_{i}, I)$.
\end{lem}
\proof
 If $PN(a_{i}, I)$ is a clique, then $a'_{i} \succ PN(a_{i}, I)$.
If $PN(a_{i}, I)$ is not a clique, let $a\in PN(a_i,I)$. From the definition of $N$, $a \succ PN(a_j,I)$ for some $a_j \in N$, which implies $a'_ja \in E(G)$. Therefore
$\hat{N} \succ PN(a_{i}, I)$ for every vertex $a_i$ of $N$.
\qed
\vskip 5 pt

\indent Now, we are ready to prove our three theorems.

\vskip 2 pt

\subsection{Proof of Theorem \ref{thm 1}}
We first consider the case when $k = \Delta + 1$. The cliques $K_k$ are isolated in $G$. Let $G_i$, $1\le i\le h$, be the components of $G$ isomorphic to $K_k$ and $G_i$, $h+1\le i\le \ell$, the other ones if any. Then $\iota _k(G)=h$ and $ir(G)\ge \ell \ge h =\iota_k$. Moreover $\iota_k(G)=ir(G)$ if and only if every component of $G$ is a clique $K_k$.


We now consider the case when $k = \Delta$.
Let $H$ be a $k$-clique of $G$.

\begin{claim}
\label{c1}
If $V(H) \cap I = \emptyset$, then $V(H) \cap S = \emptyset$
\end{claim}
\proof
We may suppose to the contrary that $V(H) \cap I = \emptyset$ but $V(H) \cap S \neq \emptyset$. Let $v \in V(H) \cap S$. Thus $deg_{H}(v) = \Delta - 1$ because $|V(H)| = k = \Delta$. Since $v \in S$, $v$ is adjacent to at least two vertices in $I$ which are not in $H$. Therefore, $deg_{G}(v) \geq deg_{H}(v) + 2 = \Delta + 1$, a contradiction.
\smallqed

\begin{claim}
\label{c2}
If $V(H) \subseteq U$ and $v \in V(H)$, then there exists $x \in N$ such that $v \succ PN(x, I)$.
\end{claim}
\proof
By Observation
1 after Theorem \ref{thm BoCo2}, $v \succ PN(x, I)$ for some $x \in B$.
Since $d_U(v)\ge d_H(v)=\Delta -1$, $|PN(x,I)|\le 1$. Thus $x\in N$
since $PN(x,I)$ is a clique.
\smallqed

\begin{claim}
\label{c4}
If $V(H) \cap (P \cup I) = \emptyset$, then, for all $v \in V(H)$, there exists $a'_{i} \in \hat{N}$, such that $va'_{i} \in E(G)$.
\end{claim}
\proof
By Claim \ref{c1}, $V(H)\cap S=\emptyset$ and thus $V(H) \subseteq U.$  By Claim \ref{c2}, $v \succ PN(a_{i}, I)$ for some $a_{i} \in N$. Thus $va'_{i} \in E(G)$.
\smallqed

\indent Let now $X = \hat{N} \cup M \cup Q \cup Z$ and let $H$ be a $k$-clique of $G$. If $V(H)\cap I \neq \emptyset$, then $H$ contains or is adjacent to a vertex of $X$ by the definition of $\hat{N}$. If $V(H)\cap P \neq \emptyset$, then $H$ is adjacent to a vertex of $X$ by Lemma \ref{l1}. If $V(H) \cap (P \cup I) = \emptyset$, then $H$ is adjacent to a vertex of $\hat{N}$ by Claim \ref{c4}. Therefore $X$ is a $K_k$-isolating set of $G$ and $\iota_{k}(G) \leq |X| = |\hat{N} \cup M \cup Q \cup Z| = ir(G)$. This completes the proof.
\qed
\vskip 10 pt


\subsection{Proof of Theorem \ref{thm 2}}
\indent Now $k=\Delta -1$. Recall that if $x\in M$, then $|PN(x,I)|\ge 2$. We let

\indent $M' =$ the vertex subset of $M$ such that $|PN(x, I)| = 2$ for all $x \in M'$,

\indent $\tilde{M} = M \setminus M'$,

\indent $S'$ = the subset of $S$ such that each vertex in $S'$ is not adjacent to any vertex in $\tilde{M} \cup Q \cup Z$.

\noindent Therefore, every vertex in $S'$ is adjacent to at least two vertices in $N \cup M'$. We let further

\indent $\mathcal{H} =$ the set of all $k$-cliques $H$ such that $V(H) \subseteq S'$.
\vskip 5 pt


\noindent To prove Theorem \ref{thm 2}, we prove in Lemma \ref{l2}
a stronger result

\noindent
\begin{lem}\label{l2}
\emph{Let $G$ be a graph with maximum degree $\Delta = k + 1$. Let  $s$ be the maximum number of $k$-cliques contained in $S'$ to which a vertex in $N \cup M'$ can be adjacent.
When $S'=\emptyset$ or does not contain any $k$-clique, let $s=0$.  Then}

$\iota_{k}(G) \leq \frac{(3\Delta - 4)ir(G)}{2\Delta - 2} - s + 1$
if $1 \leq s \leq \Delta - 2$,

$\iota _k(G)\leq ir(G)$ if $s=0$.

\end{lem}

\proof
\noindent We need the following claims.

\begin{claim}
\label{c21}
The $k$-cliques of $G[S]$ are components of $S$. In particular, the $k$-cliques in $\mathcal{H}$ are pairwise disjoint.
\end{claim}
\proof
Let $u$ be a vertex of a $k$-clique $H$ contained in $S$. Then $\Delta \ge d_G(u) \ge d_S(u)+d_I(u)\ge (k-1)+2$. Since $k=\Delta -1$,  $d_S(u)=k-1$, which shows that $H$ forms a component of $G[S]$.
\smallqed


\begin{claim}
\label{c24}
Let $H$ be a $k$-clique such that $V(H) \cap U \neq \emptyset$ and $V(H) \cap P = \emptyset$. If $u \in V(H) \cap U$, then $PN(x, I) \subseteq N_{G}(u)$ for some $x \in N \cup M'$.
\end{claim}
\proof
Since $u \in U$, there exists $x \in B$ such that $PN(x, I) \subseteq N_{G}(u)$. If $x \in \tilde{M}$, then $|PN(x, I)| \geq 3$ which implies that $deg_{P}(u) \geq 3$. Since $V(H) \cap P = \emptyset$, it follows that $deg_{G}(u) \geq deg_{H}(u) + deg_{P}(u) \geq (\Delta - 2) + 3 = \Delta + 1$, a contradiction. Therefore, $x \in N \cup M'$.
\smallqed

\begin{claim}\label{c22}
$|\mathcal{H}| \leq \frac{|N \cup M'|(\Delta - 2)}{2k}$.
\end{claim}
\proof By Claim \ref{c21}, there are exactly $2k|\mathcal{H}|$ edges from the $|\mathcal{H}|$ $k$-cliques of $S'$ to vertices in $N \cup M'$. For a vertex $x \in N \cup M'$, $deg_{I}(x) \geq 1$ because $x$ is non-isolated  in $I$ and $deg_{P}(x) \geq 1$ because $x$ has at least one private neighbor with respect to $I$. Thus, $deg_{S}(x) \leq \Delta - 2$. By a double counting, we have
\begin{align}
2k|\mathcal{H}| \leq (\Delta - 2)|N \cup M'|\notag
\end{align}
\noindent implying that $|\mathcal{H}| \leq \frac{|N \cup M'|(\Delta - 2)}{2k}$.
\smallqed

Recall that $M'$ is the vertex subset of $M$ such that $|PN(x, I)| = 2$ for all $x \in M'$. In the following, we let

\indent $M' = \{x_{1}, x_{2}, ..., x_{t}\}$ where $|M'| = t$,

\indent $PN(x_{i}, I) = \{y_{i}, w_{i}\}$,

\indent $Y = \{y_{i}: 1 \leq i \leq t\}$,

\indent $W = \{w_{i}: 1 \leq i \leq t\}$.

\noindent Moreover, we let

\indent $PN = $ the set of all private neighbors of vertices in $N$

\indent $P\tilde{M} =$ the set of all private neighbors of vertices in $\tilde{M}$,

\indent $PQ = $ the set of all private neighbors of vertices in $Q$ and

\indent $PZ = $ the set of all private neighbors of vertices in $Z$.
\vskip 5 pt

\indent When $s$, the maximum number of $k$-cliques of $ G[S']$ to which a vertex in $N \cup M'$ can be adjacent, is different from zero, let $a$ be a vertex in $N \cup M'$ which is adjacent to $s$  $k$-cliques of $\mathcal{H}$ and let $\mathcal{H}_{a}$ be the set of all cliques $H$ of $\mathcal{H}$ to which $a$ is adjacent. For each clique $H$ in $\mathcal{H} \setminus \mathcal{H}_{a}$, we choose one vertex $x_{H} \in V(H)$ and we let $X = \{a\} \cup \{x_{H} : H \in \mathcal{H} \setminus \mathcal{H}_{a}\}$. By Claim \ref{c22}, $|X| \leq \frac{|N \cup M'|(\Delta - 2)}{2k} - s + 1$. When $s=0$, we let $X=\emptyset$. Clearly,
$G[S' \setminus N[X]]$ contains no $k$-clique and thus $G[S \setminus N[X \cup \tilde{M} \cup Q \cup Z]]$ contains no $k$-clique.
\vskip 5 pt

\noindent Now, we will show that the set
\begin{center}
$T = \hat{N} \cup Y \cup X \cup \tilde{M} \cup Q \cup Z$
\end{center}
\noindent is a $K_k$-isolating set of $G$. We assume to the contrary that there exists a $k$-clique $H$ in $G - N_{G}[T]$. By the choice of $T$ and Lemma \ref{l1}, we have that
\begin{align}\label{e21}
V(H) \cap (T  \cup N \cup PN \cup M' \cup P\tilde{M} \cup PQ \cup PZ \cup (S \setminus S'))= \emptyset
\end{align}
and $H$ is not entirely contained in $S'$. Hence $V(H) \subseteq W \cup U \cup S'$ and  $V(H) \cap (W \cup U) \neq \emptyset$.

\vskip 5 pt

\indent Suppose that $V(H) \cap W \neq \emptyset$. Let $w_{i} \in V(H) \cap W$. So, $w_{i}$ is a private neighbor of a vertex $x_{i} \in M'$ whose second private neighbor $y_i$ is in $Y$. Since $Y \subset T$, $y_{i}w_{i} \notin E(G)$. By Theorem \ref{thm BoCo2}(b), $w_{i} \succ PN(x, I)$ for some $x \in I \setminus Z$. We first consider the case when $x \in N \cup M'$. So $PN(x, I) \cap (\hat{N} \cup Y) \neq \emptyset$ implying that $w_{i}$ is adjacent to a vertex in $\hat{N} \cup Y$. Thus, $w_{i} \in N_{G}[T]$, contradicting $H$ is in $G - N_{G}[T]$. Hence $x \in \tilde{M} \cup Q$ and $w_i$ has at least one neighbor $x'$ in $P\tilde M \cup PQ$. Since $V(H)\cap (P\tilde M \cup PQ)=\emptyset$, $x' \notin V(H).$
On the other hand, $x_i \in M'\subseteq B$ and by the definition of $B$, there exists a vertex $u$ in $U$ dominating $PN(x_i,I).$ Hence  $u$ is adjacent to $y_i$ and to $w_i$, and since $y_i \in T$, $u\notin V(H)$. Now $w_i$ has at least three neighbors $u, x', x_i$ not in $V(H)$ and $d_G(w_i)\ge 3+ d_H(w_i)=3+(\Delta -2)=\Delta +1$, a contradiction. Hence $V(H)\cap W=\emptyset$ and thus $V(H) \subseteq U \cup S'$ with $V(H)\cap U \neq \emptyset$.
Let $u \in V(H) \cap U$. By Claim \ref{c24}, $u$ is adjacent to a vertex in $\hat{N} \cup Y$. Therefore, $V(H) \cap N_{G}[T] \neq \emptyset$ contradicting $H$ is in $G - N_{G}[T]$. So, $T$ is a $K_{k}$-isolating set of $G$.


\noindent Hence,
\begin{align}
\iota_{k}(G) \leq |T| &\leq |\hat{N} \cup Y \cup \tilde{M} \cup Q \cup Z| + |X|\notag\\
                      &\leq ir(G) + |X|\notag\\
                      &=    ir(G) + \frac{|N \cup M'|(\Delta - 2)}{2k} - s + 1\notag\\
                      &\leq \frac{(3\Delta - 4)ir(G)}{2\Delta - 2} - s + 1\notag ~~~~ {\rm when} ~s\neq 0
\end{align}
and $\iota_k(G)\le ir(G)$ when $s=0$ since in this case, $|X|=0$.
\qed

\vskip 2 pt
From the previous lemma, $\iota _k(G)\le {\rm max} \{\frac{(3\Delta - 4)ir(G)}{2\Delta - 2}, ir(G)\}=\frac{(3\Delta - 4)ir(G)}{2\Delta - 2}$, which proves Theorem \ref{thm 2}. We see that Equality in Theorem \ref{thm 2} requires $s = 1$ as in the example given in Section 2.  In the appendix, we construct a graph satisfying the equality in Lemma \ref{l2} when $s \geq 2$.


\subsection{Proof of Theorem \ref{thm 3}}
We recall that:
\vskip 5 pt

\indent $M' =$ the subset of $M$ such that $|PN(x, I)| = 2$ for all $x \in M'$.
\vskip 5 pt

\noindent Then we define further that:
\vskip 5 pt

\indent $M'' =$ the subset of $M$ such that $|PN(x, I)| = 3$ for all $x \in M''$,
\vskip 5 pt

\indent $\overline{M} = M \setminus (M' \cup M'')$,
\vskip 5 pt

\indent $Q' =$ the subset of $Q$ such that $|PN(x, I)| = 1$ for all $x \in Q'$,
\vskip 5 pt

\indent $\tilde{Q} =$ the subset of $Q$ such that $|PN(x, I)| \geq 2$ for all $x \in \tilde{Q}$, that is $\tilde{Q} = Q \setminus Q'$,
\vskip 5 pt

\indent $S'' =$ the subset of $S$ such that every vertex in $S''$ is adjacent in $I$ to only vertices in\\
\indent $~~~~~~~$ $N \cup M' \cup M'' \cup Q'$
\vskip 5 pt

\noindent and
\vskip 5 pt

\indent $\mathcal{H}' =$ the family of $k$-cliques $H$ such that $V(H) \subseteq S''$.
\vskip 5 pt

\vskip 5 pt
The following claim is a direct consequence of Theorem \ref{thm BoCo2}(b).

\begin{claim}\label{c31}
Let $x \in M' \cup M''$ and $x' \in PN(x, I)$ such that $x'$ is not adjacent to any vertex in $PN(x, I)$. Then $x' \succ PN(y, I)$ for some $y \in A$.
\end{claim}

\vskip 5 pt

\indent When $k = \Delta - 2$ and $\Delta > 3$, it is possible that two $k$-cliques in $\mathcal{H}'$ intersect. Let $H$ and $H'$ be two different $k$-cliques in $\mathcal{H}'$ such that $V(H) \cap V(H') \neq \emptyset$. Let $s = |V(H) \cap V(H')|$ and $a \in V(H) \cap V(H')$.
Then

$d_{V(H)\cup V(H')}(a)=d_{V(H)}(a)+d_{V(H')}(a)-d_{V(H)\cap V(H')}(a)=2(k-1)-(s-1)=2k-s-1.$

\noindent
Remind that $deg_{I}(a) \geq 2$ because $a \in V(H) \subseteq S$. Thus,

$k+2=\Delta \ge d_G(a) \ge d_I(a)+d_{V(H)\cup V(H')}(a) \ge 2k-s+1$

\noindent
implying that $s\ge k-1 $.


Since $H$ and $H'$ are different, it follows that $s \leq k - 1$ which yields $s = k - 1$.
Hence two intersecting $k$-cliques of $S''$ share exactly $k-1$ vertices and each of these vertices, say $x$, has no other neighbor in $S''$ since $d_{V(H \cup H')}(x)=k+2=\Delta$. We call \emph{isolated clique} a $k$-clique of $S''$ intersecting no other one and \emph{twin} a set of two intersecting $k$-cliques of $S''$.
Note that a twin is either a $(k+1)$-clique or a $(k+1)$-clique minus one edge.
Clearly, an isolated clique has $k$ vertices while a twin has $k + 1$ vertices.
We let $t_{1}$ and $t_{2}$ be the number of isolated cliques and twins in $\mathcal{H}'$ respectively. The following claim establishes an upper bound of $t_{1} + t_{2}$.
\vskip 5 pt

\begin{claim}\label{c32}
$t_{1} + t_{2} \leq \frac{|N \cup M' \cup M'' \cup Q'|}{2}$
\end{claim}
\proof
Let $e(\mathcal{H}', N \cup M' \cup M'' \cup Q')$ be the number of edges between the vertices of the cliques in $\mathcal{H}'$ and the vertices in $N \cup M' \cup M'' \cup Q'$.
By the definition of $\mathcal{H}'$, there are at least two edges from each vertex of cliques in $\mathcal{H}'$ to $N \cup M' \cup M'' \cup Q'$. Because $\mathcal{H}'$ only consists of isolated cliques and twins, it follows that
\begin{align}\label{11}
  2k(t_{1} + t_{2}) \leq 2kt_{1} + 2(k + 1)t_{2} \leq e(\mathcal{H}', N \cup M' \cup M'' \cup Q').
\end{align}
\indent On the other hand, since every vertex $v \in N \cup M' \cup M'' \cup Q'$ is not isolated in $I$ and has at least one private neighbor in $P$, we have that $deg_{S''}(v) \leq \Delta -2.$ Hence
\begin{align}\label{12}
  e(\mathcal{H}', N \cup M' \cup M'' \cup Q') \leq (\Delta - 2)|N \cup M' \cup M'' \cup Q'| =k|N \cup M' \cup M'' \cup Q'|.
\end{align}
\noindent By Equations (\ref{11}) and (\ref{12}), $t_{1} + t_{2} \leq \frac{|N \cup M' \cup M'' \cup Q'|}{2}$ which proves Claim \ref{c32}. Note that equality in Claim \ref{c32} implies, from the first  inequality in (\ref{11}), that $t_2=0$, that is all the $k$-cliques of $S''$ are isolated.
\smallqed
\vskip 5 pt

\indent In the following, for an isolated clique or a twin $H$, we let $x_{H}$ be a vertex in $H$. We let $X = \{x_{H} : H$ is an isolated clique or a twin of $\mathcal{H'}\}$. So, $X$ is a $K_k$-isolating set of $S''$ and
\begin{align}
|X| \leq \frac{|N \cup M' \cup M'' \cup Q'|}{2}\notag
\end{align}
\noindent by Claim \ref{c32}. Recall that $\hat{N} = \{a'_{1}, a'_{2}, ..., a'_{|N|}\}$. Similarly, for a vertex $b \in M' \cup M''$, we choose one vertex $b'$ from $PN(b, I)$ and, for a vertex $c \in Q'$, we call $c'$ the unique vertex of $PN(c, I)$. Then, we let
\vskip 5 pt

\indent $\hat{M} = \{b' : b \in M' \cup M''\}$ and
\vskip 5 pt

\indent $\hat{Q} = \{c' : c \in Q'\}$.
\vskip 5 pt

\indent For the rest of this paper, we aim to prove that
\begin{center}
  $T = \hat{N} \cup \hat{M} \cup \overline{M} \cup \hat{Q} \cup \tilde{Q} \cup Z \cup X$
\end{center}
\noindent is a $K_{k}$-isolating set of $G$. We assume to the contrary that there exists a $k$-clique $H$ in $G - N_{G}[T]$. By the choice of $T$, we have
\begin{align}\label{e1}
V(H) \cap I = \emptyset.
\end{align}
\vskip 3 pt

\begin{claim}\label{c33}
$V(H) \cap U = \emptyset$.
\end{claim}
\proof
We assume to the contrary that there exists $u \in V(H) \cap U$. By Theorem \ref{thm BoCo2}(b) and the definition of $B$, $PN(x, I) \subseteq N_{G}(u)$ for some $x \in B$.
If $x \in N \cup M' \cup M''$, then $u$ is adjacent to $x'$ for some $x' \in \hat{N} \cup \hat{M} \subseteq T$ contradicting $H$ is in $G - N_{G}[T]$. Thus, $x \in \overline{M}$. If $V(H) \cap PN(x, I) \neq \emptyset$, then $H$ is adjacent to a vertex $x \in \overline{M} \subseteq T$, again contradicting $H$ is in $G - N_{G}[T]$. Thus, $V(H) \cap PN(x, I) = \emptyset$. Because $x \in \overline{M}$, we have $|PN(x, I)| \geq 4$. Thus, $deg_{PN(x, I)}(u) \geq 4$ implying that
\begin{center}
  $deg_{G}(u) \geq deg_{PN(x, I)}(u) + deg_{H}(u) \geq 4 + (k - 1) = 4 + (\Delta - 3) = \Delta + 1$,
\end{center}
\noindent a contradiction. So, $V(H) \cap U = \emptyset$ and this proves Claim \ref{c33}.
\smallqed

\begin{claim}\label{c34}
$V(H) \cap PN(b, I) = \emptyset$ for all $b \in M' \cup M''$.
\end{claim}
\proof
Assume to the contrary that there exists $b \in M' \cup M''$ such that $V(H) \cap PN(b, I) \neq \emptyset$. Let $u \in V(H) \cap PN(b, I)$. Thus, $u \in PN(b, I)$. By the definition of $M$, $u$ is adjacent to a vertex in $U$. Thus, $deg _{U}(u) \geq 1$. Clearly, $deg_{I}(u) = 1$ because $u$ is a private neighbor of $b$. By the assumption that $H$ is in $G - N_{G}[T]$, $u$ is not adjacent to $b'$ which is selected in $\hat{M}$ from $PN(b, I)$. Hence, $u$ does not dominate $PN(b, I)$. By Claim \ref{c31}, $u$ dominates PN$(y,I)$ for some $y\in A\setminus \{b\}$.
\vskip 5 pt

\indent
If $y \in (N \cup M' \cup M'') \setminus \{b\}$, then $u$ is adjacent to a vertex $y' \in \hat{N} \cup \hat{M}$ where $y'$ is selected from $PN(y, I)$, contradicting $H$ is in $G - N_{G}[T]$.

If $y \in \overline{M}$, then $V(H) \cap PN(y, I) = \emptyset$ because $\overline{M} \subseteq T.$ Since $|PN(y, I)| \geq 4$,
we have
\begin{center}
$deg_{G}(u) \geq deg_{PN(y, I)}(u) + deg_{H}(u) \geq 4 + (k - 1) = \Delta + 1$,
\end{center}
\noindent a contradiction.

Finally, we consider the case when $y \in Q$. Because $\hat{Q} \subseteq T$, it follows that $y \in \tilde{Q}$. In this case $|PN(y, I)| \geq 2$. Similarly, $V(H) \cap PN(y, I) = \emptyset$ because $\tilde{Q} \subseteq T$. By Claim \ref{c33} and Equation (\ref{e1}),
\begin{center}
$deg_{G}(u) \geq deg_{PN(y, I)}(u) + deg_{H}(u) + deg_{U}(u) + deg_{I}(u) \geq 2 + (k - 1) + 1 + 1 = \Delta + 1$,
\end{center}
\noindent a last contradiction
proving Claim \ref{c34}.
\smallqed

\indent Now, for a $k$-clique $H$ which is in $G - N_{G}[T]$, we have by Equation (\ref{e1}), Claim \ref{c34} and Lemma \ref{l1} that
\begin{center}
$V(H) \cap (I \cup P) = \emptyset$.
\end{center}
\noindent Moreover, by Claim \ref{c33}, $V(H) \subseteq S$. Because $\overline{M} \cup \tilde{Q} \cup Z \subseteq T$, it follows that every vertex in $H$ is adjacent to only vertices in $N \cup M' \cup M'' \cup Q'$. Therefore, $V(H) \subseteq S''$, namely, $H \in \mathcal{H}'$. It follows that $H$ contains or is adjacent to a vertex in $X$, a contradiction since $X \subset T$. Therefore, $T$ is an $K_{k}$-isolating set of $G$.
\vskip 5 pt

\noindent By the minimality of $\iota_{k}(G)$ and Claim \ref{c32}, we have
\begin{align}
\iota_{k}(G) \leq |T| &=    |\hat{N} \cup \hat{M} \cup \overline{M} \cup \hat{Q} \cup \tilde{Q} \cup Z \cup X|\notag\\
                      &=    |\hat{N} \cup \hat{M} \cup \overline{M} \cup \hat{Q} \cup \tilde{Q} \cup Z| + |X|\notag\\
                      &\leq ir(G) + \frac{|N \cup M' \cup M'' \cup Q'|}{2}\notag\\
                      &\leq ir(G) + \frac{ir(G)}{2} = \frac{3ir(G)}{2}.\notag
\end{align}
\noindent This completes the proof.
\qed

\noindent \textbf{Appendix}
\vskip 5 pt

\indent We give a construction of graphs satisfying the bound of Lemma \ref{l2} when $s \geq 2$ since an example with $s = 1$ was already shown in Section \ref{sec2}. First, we  construct the graph $S(k)$ which is obtained from two paths $x_{1}x_{2}x_{3}x_{4}$ and $x_{5}x_{6}x_{7}x_{8}$ of length 3 and five $k$-cliques $K^{0}_{k}, K^{1}_{k}, ..., K^{4}_{k}$ by joining edges as follows:
\begin{itemize}
  \item join $x_{1}, x_{4}, x_{5}$ and $x_{8}$ to every vertex of $K^{1}_{k}, K^{2}_{k}, K^{3}_{k}$ and $K^{4}_{k}$, respectively.
\end{itemize}
\noindent Further, we let $V(K^{0}_{k}) = \{z_{1}, ..., z_{k}\}$. Then,
\begin{itemize}
  \item join $x_{2}$ to $z_{1}, ..., z_{k - 1}$, $x_{3}$ to $z_{2}, ..., z_{k}$ and join $x_{6}$ to $z_{1}$ and $z_{k}$.
\end{itemize}

\begin{figure}[H]
\centering
\definecolor{ududff}{rgb}{0.30196078431372547,0.30196078431372547,1}
\resizebox{0.7\textwidth}{!}{%
\begin{tikzpicture}[line cap=round,line join=round,>=triangle 45,x=1cm,y=1cm]

\draw [line width=0.5pt] (3,7) circle (0.4249952709310901cm);
\draw [line width=0.5pt] (5,7) circle (0.4165789814669837cm);
\draw [line width=0.5pt] (8,7) circle (0.4117977543525712cm);
\draw [line width=0.5pt] (10,7) circle (0.4351550672108444cm);
\draw [line width=0.5pt] (3,4)-- (5,4);
\draw [line width=0.5pt] (8,4)-- (10,4);
\draw [line width=0.5pt] (3,5)-- (3,4);
\draw [line width=0.5pt] (5,5)-- (5,4);
\draw [line width=0.5pt] (8,5)-- (8,4);
\draw [line width=0.5pt] (10,5)-- (10,4);
\draw [line width=0.5pt] (3,4)-- (3,2);
\draw [line width=0.5pt] (3,4)-- (4,2);
\draw [line width=0.5pt] (3,4)-- (5,2);
\draw [line width=0.5pt] (3,4)-- (6,2);
\draw [line width=0.5pt] (3,4)-- (7,2);
\draw [line width=0.5pt] (3,4)-- (8,2);
\draw [line width=0.5pt] (3,4)-- (9,2);
\draw [line width=0.5pt] (5,4)-- (10,2);
\draw [line width=0.5pt] (5,4)-- (9,2);
\draw [line width=0.5pt] (5,4)-- (8,2);
\draw [line width=0.5pt] (5,4)-- (7,2);
\draw [line width=0.5pt] (5,4)-- (6,2);
\draw [line width=0.5pt] (5,4)-- (5,2);
\draw [line width=0.5pt] (5,4)-- (4,2);
\draw [line width=0.5pt] (8,4)-- (10,2);
\draw [line width=0.5pt] (8,4)-- (3,2);
\draw [rotate around={0.1613961351203842:(6.515711205836903,2.0107937378855394)},line width=0.5pt] (6.515711205836903,2.0107937378855394) ellipse (4.0758857740632735cm and 0.7779177632633627cm);
\draw [line width=0.5pt] (3,5)-- (2.625279801154102,6.7994862775501845);
\draw [line width=0.5pt] (3,5)-- (3.3791912048897044,6.808075508472512);
\draw [line width=0.5pt] (5,5)-- (4.62496106726536,6.818660961915721);
\draw [line width=0.5pt] (5,5)-- (5.368082569766187,6.80492752722542);
\draw [line width=0.5pt] (8,5)-- (7.6269467616973445,6.825619175705232);
\draw [line width=0.5pt] (8,5)-- (8.365693150162242,6.810669837547641);
\draw [line width=0.5pt] (10,5)-- (9.609010192305838,6.8089845482727736);
\draw [line width=0.5pt] (10,5)-- (10.392464335651404,6.8120327747691745);
\begin{scriptsize}
\draw [fill=black] (3,2) circle (2pt);
\draw [fill=black] (4,2) circle (2pt);
\draw [fill=black] (5,2) circle (2pt);
\draw [fill=black] (6,2) circle (2pt);
\draw [fill=black] (7,2) circle (2pt);
\draw [fill=black] (8,2) circle (2pt);
\draw [fill=black] (9,2) circle (2pt);
\draw [fill=black] (10,2) circle (2pt);
\draw [fill=black] (3,4) circle (2pt);
\draw[color=black] (3.2,4.2) node {$x_{2}$};
\draw [fill=black] (5,4) circle (2pt);
\draw[color=black] (5.2,4.2) node {$x_{3}$};
\draw [fill=black] (8,4) circle (2pt);
\draw[color=black] (8.2,4.2) node {$x_{6}$};
\draw [fill=black] (10,4) circle (2pt);
\draw[color=black] (10.2,4.2) node {$x_{7}$};
\draw [fill=black] (3,5) circle (2pt);
\draw[color=black] (3.3,5.2) node {$x_{1}$};
\draw [fill=black] (5,5) circle (2pt);
\draw[color=black] (5.3,5.2) node {$x_{4}$};
\draw [fill=black] (8,5) circle (2pt);
\draw[color=black] (8.3,5.2) node {$x_{5}$};
\draw [fill=black] (10,5) circle (2pt);
\draw[color=black] (10.3,5.2) node {$x_{8}$};

\draw[color=black] (10.7,7) node {$K^{4}_{8}$};
\draw[color=black] (8.7,7) node {$K^{3}_{8}$};
\draw[color=black] (5.7,7) node {$K^{2}_{8}$};
\draw[color=black] (3.7,7) node {$K^{1}_{8}$};

\end{scriptsize}
\end{tikzpicture}
 }%
\vskip -0.25 cm
\caption{The graph $S(8)$.}
\label{s8}
\end{figure}
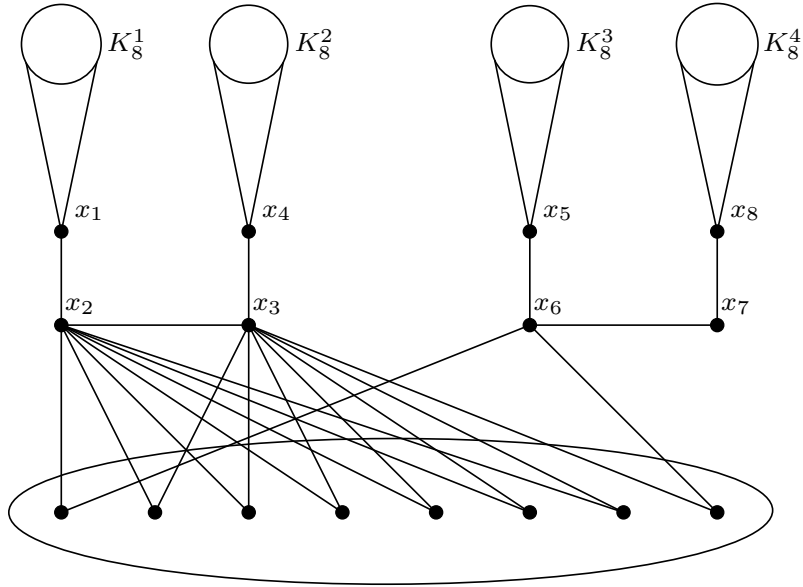

\noindent Note that $deg_{S(k)}(x_{2}) = deg_{S(k)}(x_{3}) = k + 1 = \Delta$ while $d_{S(k)}(x_{6}) = 4$ and $d_{S(k)}(x_{7}) = 2$. An example of the graph $S(k)$ when $k = 8$ is illustrated by Figure \ref{s8}.

\begin{figure}[H]
\centering
\definecolor{ududff}{rgb}{0.30196078431372547,0.30196078431372547,1}
\resizebox{0.7\textwidth}{!}{%
\begin{tikzpicture}[line cap=round,line join=round,>=triangle 45,x=1cm,y=1cm]

\draw [rotate around={-0.044968085642235096:(1.8043864578692221,2.001095335166135)},line width=1pt] (1.8043864578692221,2.001095335166135) ellipse (1.5768080735113994cm and 0.7338837387166961cm);
\draw [rotate around={0:(5.102151781791028,2.0021906703322707)},line width=1pt] (5.102151781791028,2.0021906703322707) ellipse (0.8686627090908359cm and 0.5099498435457904cm);
\draw [rotate around={0.08967734144480548:(8.699820377295834,2.001095335166135)},line width=1pt] (8.699820377295834,2.001095335166135) ellipse (0.8726383110619215cm and 0.5212944098063566cm);
\draw [line width=0.5pt] (1,4)-- (0.408772915738447,2.0021906703322707);
\draw [line width=0.5pt] (1,4)-- (4.398926313514559,2.0021906703322707);
\draw [line width=0.5pt] (1,4)-- (8,2);
\draw [line width=0.5pt] (1,4)-- (0.8,2);
\draw [line width=0.5pt] (3,4)-- (0.8,2);
\draw [line width=0.5pt] (3,4)-- (1.2,2);
\draw [line width=0.5pt] (3,4)-- (1.6,2);
\draw [line width=0.5pt] (3,4)-- (2.8,2);
\draw [line width=0.5pt] (3,4)-- (2.4,2);
\draw [line width=0.5pt] (3,4)-- (2,2);
\draw [line width=0.5pt] (1,4)-- (3,4);
\draw [line width=0.5pt] (1,4)-- (1,5);
\draw [line width=0.5pt] (3,4)-- (3,5);
\draw [line width=0.5pt] (1,7) circle (0.4239758018579873cm);
\draw [line width=0.5pt] (3,7) circle (0.42176406919404186cm);
\draw [line width=0.5pt] (1,5)-- (0.6893622616798271,6.711452471692355);
\draw [line width=0.5pt] (1,5)-- (1.3111648700833036,6.7120210004407195);
\draw [line width=0.5pt] (3,5)-- (2.6839461253328083,6.720724368461959);
\draw [line width=0.5pt] (3,5)-- (3.3098632000649517,6.71387113513554);
\draw [line width=0.5pt] (3,4)-- (3.2,2);
\begin{scriptsize}
\draw [fill=black] (1,4) circle (2pt);
\draw[color=black] (1.2,4.2) node {$a_{2}$};
\draw [fill=black] (3,4) circle (2pt);
\draw[color=black] (3.2,4.2) node {$a_{3}$};
\draw [fill=black] (0.408772915738447,2.0021906703322707) circle (2pt);
\draw [fill=black] (3.2,2) circle (2pt);
\draw [fill=black] (4.398926313514559,2.0021906703322707) circle (2pt);
\draw [fill=black] (8,2) circle (2pt);
\draw [fill=black] (0.8,2) circle (2pt);
\draw [fill=black] (1.2,2) circle (2pt);
\draw [fill=black] (1.6,2) circle (2pt);
\draw [fill=black] (2,2) circle (2pt);
\draw [fill=black] (2.4,2) circle (2pt);
\draw [fill=black] (2.8,2) circle (2pt);
\draw [fill=black] (1,5) circle (2pt);
\draw[color=black] (1.2,5.2) node {$a_{1}$};
\draw [fill=black] (3,5) circle (2pt);
\draw[color=black] (3.2,5.2) node {$a_{4}$};
\draw[color=black] (7,2) node {$...$};

\draw[color=black] (2,1) node {$\tilde{K}^{1}_{8}$};
\draw[color=black] (9,1) node {$\tilde{K}^{6}_{8}$};
\draw[color=black] (5,1) node {$\tilde{K}^{2}_{8}$};

\draw[color=black] (1.7,7) node {$K'_{8}$};
\draw[color=black] (3.7,7) node {$K''_{8}$};

\end{scriptsize}
\end{tikzpicture}
 }%
\vskip 1 cm
\caption{The graph $F(8, 6)$.}
\label{f86}
\end{figure}
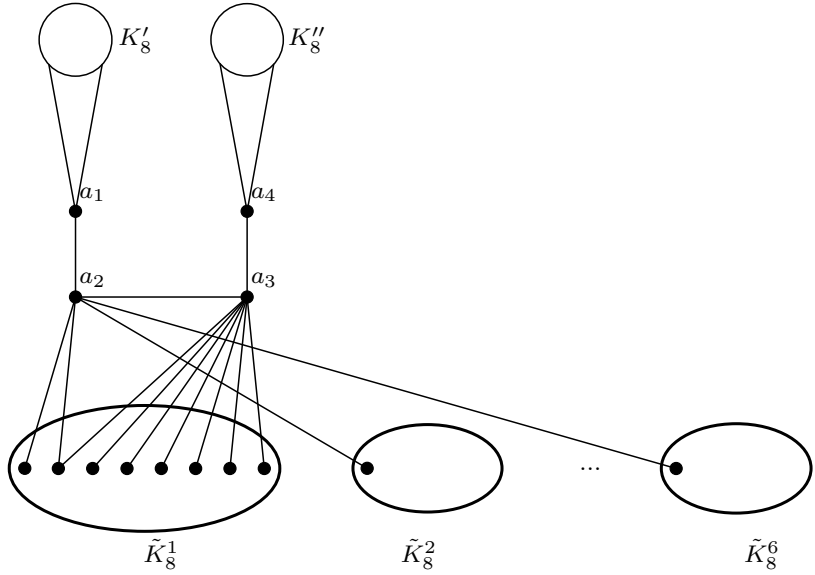

\noindent Further, we construct the graph $F(k, s)$ which is obtained from a path $a_{1}a_{2}a_{3}a_{4}$ of length $3$ and $s + 2$ of $k$-cliques $K'_{k}, K''_{k}$ and $\tilde{K}^{1}_{k}, ..., \tilde{K}^{s}_{k}$ by putting edges as follows:
\begin{itemize}
  \item join $a_{1}$ and $a_{4}$ to every vertex of $K'_{k}$ and $K''_{k}$, respectively,
  \item join $a_{2}$ to two vertices of $\tilde{K}^{1}_{k}$ and to one vertex of each $\tilde{K}^{2}_{k}, ..., \tilde{K}^{s}_{k}$ and
  \item join $a_{3}$ to $\Delta - 2$ vertices in $\tilde{K}^{1}_{k}$ to share exactly one common neighbor with $a_{2}$.
\end{itemize}
\noindent It is worth noting that exactly one vertex in $\tilde{K}^{1}_{k}$ has degree $\Delta$ and the other $k - 1$ vertices have degree $\Delta - 1$. An example the graph $F(k, s)$ when $k = 8$ and $s = 6$ is illustrated by Figure \ref{f86}.

\indent Now, we are ready to construct the graph $H(k, s, l)$ from one copy of $F(k, s)$ and $l$ copies $S_{1}, ..., S_{l}$ of $S(k)$ where $l$ is given by the equation
\begin{align}\label{5}
(s-1)(2k-1) =  (2\Delta - 6)(l - 1) + r
\end{align}
\noindent for some $0 \leq r < 2\Delta -6$. We may relabel $x_{j}$ in $S(k)$ to be $x^{i}_{j}$ in $S_{i}$ for all $1 \leq j \leq 8$ and $1 \leq i \leq l$. Further, we let
\vskip 5 pt

\noindent $L_1=\bigcup_{i=2}^s\tilde{K}_k^i$ in the copy of $F(k,s)$ and
\vskip 5 pt

\noindent $L_{2} = \bigcup^{l-1}_{i= 1}\{x^{i}_{6}, x^{i}_{7}\} \cup\{x_6^l\}$ in the $l$ copies of $S(k)$.
\vskip 5 pt

\noindent Then we add edges as follows, in such a way that the maximum degree of $H(k,s,l)$ remains equal to $\Delta =k+1$.
\begin{itemize}
  \item join $x^{l}_{7}$ to all the $k - 1$ vertices of degree $\Delta - 1$ in $V(\tilde{K}^{1}_{k})$,
  \item join $x^{l}_{6}$ to $r'=min \{r,\Delta -4\}$ vertices in $L_1$,
  \item for $1\le i\le l-1$, join $x_6^i$ to $\Delta -4$ vertices of  $L_1$ and $x_7^i$ to $\Delta - 2$ vertices of $L_1$ in such a way that every vertex in $L_{1}$ is adjacent to at most two vertices in $L_{2}\cup \{a_2,a_3\}$ to respect the value of the maximum degree.
\end{itemize}
\noindent The adjacency condition above is possible since $s-1+(2\Delta -6)(l-1)+r' \le 2(s-1)k$ from Equation (\ref{5}).

\vskip 5 pt





We show that the graph $H(k, s, l)$ when $k = 8, s = 6$ and $l = 7$ satisfies the upper bound of Lemma \ref{l2}. For the sake of convenience, we let $H = H(8, 6, 7)$. Clearly, $\{a_{2}, x^{1}_{2}, ..., x^{7}_{2}\} \cup \{a_{1}, a_{4}\} \cup (\cup^{7}_{i = 1}\{x^{i}_{1}, x^{i}_{4}, x^{i}_{5}, x^{i}_{8}\})$ is a $K_8$-isolating set of $H$. Thus $\iota_{8}(H) \leq 38$. Let $S$ be a smallest $\iota_{8}$-set of $G$. To be adjacent to all the $8$-cliques, $S$ intersects each disjoint subset as follows: $(i)$ $S \cap (\{a_{1}\} \cup V(K'_{8})) \neq \emptyset, S \cap (\{a_{4}\} \cup V(K''_{8})) \neq \emptyset$, $(ii)$ $S \cap (\{x^{i}_{j}\} \cup V(K^{i, \lfloor\frac{j + 1}{2}\rfloor)}_{8})$ for all $1 \leq i \leq 7$ and $j \in \{1, 4, 5, 8\}$, $(iii)$ $S \cap (\{a_{2}, a_{3}, x^{7}_{7}\} \cup V(\tilde{K}^{1}_{8})) \neq \emptyset$ and $(iv)$ $S \cap (\{x^{i}_{2}, x^{i}_{3}, x^{i}_{6}\} \cup V(K^{i, 0}_{k})) \neq \emptyset$ for all $1 \leq i \leq 7$. By $(i)$ - $(iv)$, $S$ has at least $2 + 28 + 1 + 7 = 38$ vertices implying that $\iota_{8}(H) = 38$.
\vskip 5 pt

\indent We next let $I$ be a smallest maximal irredundant set of $H$. Hence, every vertex in $I$ has a private neighbor with respect to $I$. If $I \cap (\{a_{1}, a_{2}\} \cup V(K'_{8})) = \emptyset$, then $a_{1}$ is not a private neighbor of any vertex in $I$. Thus, for a vertex $h \in V(K'_{8})$, $a_{1}$ is a private neighbor of $h$ with respect to $I \cup \{h\}$. Therefore, $I \cup \{h\}$ is an irredundant set of $H$ containing $I$ contradicting maximality of $I$. So, $I \cap (\{a_{1}, a_{2}\} \cup V(K'_{8})) \neq \emptyset$. Similarly, $I \cap (\{a_{3}, a_{4}\} \cup V(K''_{8})) \neq \emptyset, I \cap (\{x^{i}_{j}, x^{i}_{j + 1}\} \cup V(K^{i, \lfloor\frac{j + 1}{2}\rfloor}_{8})) \neq \emptyset$ and $I \cap (\{x^{i}_{j'}, x^{i}_{j' - 1}\} \cup V(K^{i, \lfloor\frac{j' + 1}{2}\rfloor}_{8})) \neq \emptyset$ for all $j \in \{1, 5\}, j' \in \{4, 8\}$ and $1 \leq i \leq 7$. This implies that $ir(H) = |I| \geq 1 + 1 + 28 = 30$. On the other hand, we observe that $\{a_{2}, a_{3}\} \cup (\cup^{7}_{i = 1}\{x^{i}_{2}, x^{i}_{3}, x^{i}_{6}, x^{i}_{7}\})$ is an irredundant set of $G$. By minimality of $ir(H)$, we have $ir(H) \leq 30$ implying that $ir(H) = 30$. We see that, in the graph $H = H(8, 6, 7)$, we have $k = 8, \Delta = 9, s = 6$ and $l = 7$. Thus, $\iota_{8}(H) = 38 = 30 + \lfloor\frac{30(9 - 2)}{2(8)}\rfloor - 6 + 1 = ir(H) + \lfloor\frac{ir(H)(\Delta - 2)}{2k}\rfloor - s + 1$ satisfying the bound of Lemma \ref{l2} when $1 \leq s \leq \Delta - 2$.

\end{document}